\documentclass[11pt]{article}
\usepackage{amsfonts}
\usepackage{epsf}

\newcommand{\be}{\begin{equation}}
\newcommand{\ee}{\end{equation}}
\newcommand{\bea}{\begin{eqnarray}}
\newcommand{\eea}{\end{eqnarray}}
\newcommand{\beann}{\begin{eqnarray*}}
\newcommand{\eeann}{\end{eqnarray*}}
\newcommand{\beaa}{\begin{eqnarray*}}
\newcommand{\eeaa}{\end{eqnarray*}}

\newcommand{\wphi}{\widehat{\phi}}
\newcommand{\wpsi}{\widehat{\psi}}
\newcommand{\wu}{\widehat{u}}
\newcommand{\wv}{\widehat{v}}
\newcommand{\ww}{\widehat{w}}
\newcommand{\wf}{\widehat{f}}
\newcommand{\wg}{\widehat{g}}

\newcommand{\eq}[1]{(\ref{#1})}

\newcommand{\Rl}{\mathbb{R}}
\newcommand{\Cx}{\mathbb{C}}
\newcommand{\Nl}{\mathbb{N}}
\newcommand{\Ir}{\mathbb{Z}}
\newcommand{\Ts}{\mathbb{T}}

\newcommand{\fd}[2]{\frac{\delta #1}{\delta #2}}
\newcommand{\sgn}{\mathrm{sgn}\,}

\newcommand{\scaling}{\eta}
\newcommand{\hilbert}{\mathbf{H}}
\newcommand{\hilb}{\mathbf{H}}

\newcommand{\ham}{\mathcal{H}}
\newcommand{\mom}{\mathcal{P}}

\renewcommand{\phi}{\varphi}

\newtheorem{proposition}{Proposition}
\newtheorem{theorem}{Theorem}

\begin{document}

\title{Short-Time Existence for\\
Scale-Invariant Hamiltonian Waves}

\author{John K. Hunter\footnote{Partially supported by the NSF under grant numbers DMS-0309648 and FRG-0243622}
\\ Department of Mathematics  \\
University of California at Davis\\
Davis, CA 95616, USA}

\date{June 28, 2005}

\maketitle

\begin{abstract}
We prove short-time existence of smooth solutions for a class of nonlinear,
and in general spatially nonlocal, Hamiltonian evolution equations that describe the self-interaction of weakly nonlinear
scale-invariant waves. These equations include ones that describe
weakly nonlinear hyperbolic surface waves, such as nonlinear Rayleigh waves
in elasticity.
\end{abstract}

\section{Introduction}
\setcounter{equation}{0}

In this paper we prove the short-time existence of solutions of
a class of quadratically nonlinear evolution equations in one space-dimension. The equations
are, in general, spatially nonlocal; an example is the following singular integro-differential
equation for $\phi(x,t)$
\be
\phi_t  +\frac{1}{2}
\hilbert\left[\psi^2\right]_{xx}
+ \psi \phi_{xx}
 = 0,
 \label{surf_ex}
\ee
where $\hilbert$ denotes the Hilbert transform with respect
to $x$ and $\psi = \hilbert[\phi]$. According to the result
we prove here, the initial value problem for this equation
has a unique short-time solution in the Sobolev space $H^s$
when $s > 5/2$.

This class of equations arises from a Hamiltonian description \cite{zakharov} of
weakly nonlinear, translation-invariant, scale-invariant waves \cite{parker}.
By `scale-invariant' we mean that the wave motion is invariant under
space-time rescalings $x\mapsto \scaling x$, $t\mapsto\scaling t$,
for all $\scaling > 0$, with the dependent field variables scaling according to
their appropriate space-time dimensions. As a result, the wave motion
does not depend on any intrinsic length or time scales, and any
space-time parameters on which it does depend
must have the dimension of velocity.\footnote{Wave motions,
such as gravity water waves on deep
water, that depend only on an acceleration rather than a velocity are not scale-invariant
in the sense we use the term here.} One consequence of this scale-invariance
is that the wave motion is non-dispersive.

A specific physical motivation for studying this class of equations
comes from quasi-linear hyperbolic surface waves,
such as Rayleigh waves in elasticity.\footnote{The term `surface wave' is also used to refer to water waves,
whose properties are different from the non-dispersive surface waves we discuss here.}
Surface waves that satisfy half-space or planar interface problems for systems of conservation
laws (without source terms) are scale-invariant because neither the PDE nor the geometry define an intrinsic
length scale.

We refer to surface waves, like Rayleigh waves,
that decay exponentially away from the boundary on which they propagate,
as `genuine' surface waves, in contrast
to `radiative' or `leaky' surface waves that propagate bulk waves into the interior.
Weakly nonlinear `genuine' surface waves are described asymptotically
by nonlocal equations, of which \eq{surf_ex} is an example,
whereas purely `radiative' surface waves are typically described by standard local PDEs
(see \cite{artola} for an example).
One of our main goals is to prove short-time existence for such nonlocal asymptotic surface wave equations, and
it is useful to put them in the context of a broader Hamiltonian framework.

These surface wave equations are of interest, in particular, because they describe the formation of singularities
in quasi-linear hyperbolic waves on the boundary of a half-space, a mechanism that differs from the familiar
one of shock formation in the interior of a domain. Their study should also shed light on the
nonlinear well-posedness of IBVPs for hyperbolic PDEs, and on the stability of shocks and other interfaces,
when the uniform Lopatinski condition fails as a result of the presence of `genuine' surface waves.
(See \cite{ali}, \cite{serre} for further discussion, and \cite{sable-tougeron} for a local existence result for
nonlinear boundary value problems in elasticity that permit Rayleigh waves.)

The canonical spectral form of the wave equation we study here
is given in \eq{specee}. This equation describes the evolution of a wavefield subject to
three-wave resonant interactions between wavenumbers $\{k, m, n\}$
such that $k+m+n=0$. These three-wave resonances always occur for scale-invariant waves, since
they are non-dispersive. Equation \eq{specee} depends on a given kernel, or interaction coefficient,
$T : \Ir^3 \to \Cx$.
Provided that $T(k,m,n)$ is not identically zero when $k+m+n=0$, we expect
that the dominant weakly nonlinear effects are those given in \eq{specee}; thus,
this equation provides a general description of nonlinear, scale-invariant waves in one space
dimension when higher-order resonances can be neglected.

The Hamiltonian formalism used to derive \eq{specee} does not imply that the resulting equations
are well-posed. Short-time existence is not obvious, and suitable
bounds on the interaction coefficients appear to be required. (See
\cite{hunter} for some examples of related nonlocal, but non-Hamiltonian,
evolution equations that are ill-posed, and the discussion of compacton-type equations in
Section~\ref{sec:compacton}.)
Consequently, the result we prove here provides a useful validation of the corresponding equations
as consistent models for the self-interaction of nonlinear waves.

Solutions of these equations typically form singularities in finite time, so that, in general, a short-time
existence result for smooth solutions is the best one can hope for.
Singularity formation and the global existence of appropriate weak solutions
are interesting questions, but we do not address them here. Unlike local existence,
global existence presumably depends on the specific kernel,
making a general global existence theory unlikely.

In proving our existence result, we find it convenient
to study a noncanonical form of \eq{specee}, given in
\eq{ncspatee}, which may be regarded as a nonlocal
generalization of the inviscid Burgers equation.
The proof uses a standard Galerkin procedure, although
we work in Fourier space because of the spatially nonlocal structure of the evolution equation.

The Hamiltonian structure of \eq{specee} imposes
certain reality and symmetry conditions on the kernel $T$, given in \eq{propT}.
In particular, it implies that $T$ is symmetric under under cyclic
permutations of its arguments, and this property plays a crucial role in our
short-time existence proof.

Scale-invariance further implies that $T$ is a homogeneous function whose
degree $\nu$  is a parameter that characterizes the space-time scaling invariance of the equation.
This parameter takes the value $3/2$ for bulk waves and $2$ for surface waves, but our analysis
applies to kernels of arbitrary degree.

It is convenient to treat the cases $\nu \ge 3/2$ and $\nu < 3/2$
separately. For $\nu \ge 3/2$,
we assume that the interaction
coefficients satisfy the bound in \eq{estimate}. This condition
appears to be necessary to develop a short-time $H^s$-existence theory
in this case, and it is satisfied
by the kernels that arise in surface wave problems.
Our existence result is stated in Theorem~\ref{th:locex}, and the
main estimate is given in Proposition~\ref{prop:estimate}.
One quantity of interest is the number of square-integrable
spatial derivatives required to prove short-time existence.
As we show, short-time existence holds in $H^s$ for the noncanonical
equation \eq{ncspatee} when $s > \nu$.

An existence result for $\nu < 3/2$ under the
assumption that the interaction coefficients satisfy the bound in \eq{negestimate} ---
a weaker condition than \eq{estimate} --- is given in Theorem~\ref{th:locexneg}.

In Section~\ref{sec:ee}, we write out the general equation in various forms (canonical and noncanonical,
spatial and spectral). In Section~\ref{sec:ex}, we give some examples that arise from specific choices
of the kernel $T$. These include the inviscid
Burgers equation, the Hunter-Saxton equation, asymptotic equations for nonlinear surface waves,
and equations related to compacton equations. In Section~\ref{sec:locex}, we prove the local existence result,
and in Section~\ref{sec:ineq} we state an inequality used in the proof.

\section{The evolution equation}
\label{sec:ee}

We consider an equation for a scalar-valued function $\phi(x,t)$ of the form
\be
\phi_t + \hilb\left[b(\phi,\phi)\right] = 0.
\label{eeq}
\ee
Here, $\hilb$ denotes the Hilbert transform with respect to the one-dimensional
spatial variable $x$, and $b(\phi,\phi)$ is a quadratically nonlinear function.
The skew-adjoint Hilbert transform plays the role of a Hamiltonian operator,
and $b(\phi,\phi)$ is the variational derivative
of a cubically nonlinear Hamiltonian functional.

The specific structure of the equation is most easily described in Fourier space.
We study spatially periodic solutions for definiteness, so that
\[
\phi : \Ts \times I \to \Rl
\]
where $\Ts$ is the circle of length $2\pi$, and $I =(-t_\ast,t_\ast)$ is a time interval.
Most of the results extend in a straightforward
way to equations defined on $\Rl$ instead of $\Ts$.

We denote the Fourier coefficients of $\phi$
by $\wphi : \Ir\times I \to \Cx$, where
\beann
&&\phi(x,t) = \sum_{k=-\infty}^{\infty} \wphi(k,t) e^{ikx},\\
&&\wphi(k,t) = \frac{1}{2\pi}\int_{\Ts} \phi(x,t) e^{-ikx}\,dx.
\eeann
The Hilbert transform is defined by
\beann
\hilbert[\phi](x,t) &=& - \frac{1}{2\pi} \mathrm{p.v.} \int_{\Ts} \cot\left(\frac{x-y}{2}\right) \phi(y,t)\,dy,\\
\widehat{\hilbert[\phi]}(k,t) &=& i\sgn k\, \wphi(k,t),
\eeann
where
\[
\sgn k = \left\{\begin{array}{rl}
1 & \mbox{$k > 0$},
\\
0 & \mbox{$k = 0$},
\\
-1 & \mbox{$k < 0$}.
\end{array}\right.
\]

The bilinear form $b$ is defined in Fourier space by
\be
\widehat{b(\phi,\psi)}(k,t) = \sum_{n=-\infty}^\infty T(-k, k - n, n) \wphi(k-n,t) \wpsi(n,t),
\label{defF}
\ee
where $T : \Ir^3 \to \Cx$ is a given kernel.
Only the values of $T(k,m,n)$ with
$
k+m+n=0
$
appear in \eq{defF},
so we could omit one of the arguments,
but it is convenient to include all of them to exhibit the Hamiltonian symmetry property
stated below. It then follows from \eq{eeq} and \eq{defF} that $\wphi(k,t)$ satisfies the
evolution equation
\be
\wphi_t(k,t) + i\sgn k \sum_{n=-\infty}^\infty T(-k, k - n, n) \wphi(k-n,t) \wphi(n,t) = 0.
\label{specee}
\ee

We assume that $T : \Ir^3\to \Cx$
satisfies the following reality, symmetry, Hamiltonian, and scaling conditions:
\bea
&&T(k, m , n)  =T^\ast(-k, -m , -n),
\nonumber\\
&&T(k, m , n) = T(k, n , m),
\label{propT}\\
&&T(k, m , n) = T(m , n, k),
\nonumber
\\
&&T(\scaling k, \scaling m , \scaling n) = \scaling^\nu T(k, m , n)\qquad \mbox{for all $\scaling >0$}.
\nonumber
\eea
Here, a star denotes the complex-conjugate. The symmetry of $T(k,m,n)$ under cyclic
permutations of $\{k,m,n\}$ gives \eq{specee} a Hamiltonian structure.
For simplicity, we also assume that
\[
T(0,m,n) = T(k,0,n) = T(k,m,0) = 0,
\]
so that there are no mean-field interactions.

The homogeneity of $T$ of degree $\nu$ corresponds to
an invariance of \eq{eeq} under rescalings
\be
t\mapsto \scaling t,\quad x\mapsto \scaling x,
\quad\phi\mapsto \scaling^{\nu-1} \phi.
\label{phi_scale}
\ee
Dimensional analysis \cite{parker} shows that the parameter $\nu$ is related to the co-dimension
$(n-d)$ of the waves by
\be
\nu = \frac{n-d+3}{2}.
\label{defnu}
\ee
Here,  $d$ is the dimension
of the wavenumber space, and $n$ is the spatial dimension of the parameter $\rho$ that sets the mass dimension
of the wave motion; for example, if this parameter is a density, then $n$ is the dimension
of the space in which the waves propagate.
It follows that $\nu= 3/2$ for bulk waves with $d=n$, and $\nu = 2$ for surface waves
with $d=n-1$.

The dimensional analysis also shows that if $\phi$ is a canonical variable, then
$\phi/\rho$ has dimension $L^{\nu-1/2}T^{-1/2}$, where $L$ and $T$ denote space and time dimensions, respectively.
This explains the scaling of $\phi$ by $\scaling^{\nu-1}$ in \eq{phi_scale}.

\subsection{Hamiltonian structure}

For simplicity, we consider functions whose spatial mean is zero, so that $\wphi(0,t) = 0$.
(This involves no loss of generality, since $\wphi(0,t)$ is constant in time.) We define
the Hamiltonian
\beann
\ham(\wphi,\wphi^\ast) &=& \sum_{m,n=1}^\infty \left\{T(-m-n, m, n)
\wphi^\ast(m+n,t)\wphi(m,t) \wphi(n,t)\right. \\
&&\qquad\qquad+
\left. T^\ast(-m-n, m, n)\wphi^\ast(m,t) \wphi^\ast(n,t) \wphi(m+n,t) \right\}.
\eeann
The complex canonical form of Hamilton's equations,
\[
i \wphi_t(k,t) = \fd{\ham}{\wphi^*}(k,t)\qquad\mbox{for $k \in \Nl$},
\]
then gives
\beann
&&i\wphi_t(k,t) = \sum_{n=1}^{k-1} T(-k, k-n, n) \wphi(k-n,t) \wphi(n,t)\\
&&\qquad\qquad+ 2\sum_{n=1}^{\infty} T^\ast(-k-n, k, n) \wphi(k+n,t) \wphi^\ast(n,t).
\eeann
The first term on the right-hand side of this equation describes the upward transfer of
energy by nonlinear interactions to higher wavenumbers, while the second term describes the reverse downward transfer.
These equations may be rewritten as \eq{specee} for $k \in \Ir$
by use of the properties of $T$ and the fact that $\wphi^\ast(k,t) = \wphi(-k,t)$.
Thus, equation \eq{specee} is Hamiltonian and
$
\left\{\wphi(k,t), \wphi^\ast(k,t) \mid k\in\Nl\right\}
$
are complex canonical variables.

The Hamiltonian is conserved in time on smooth solutions,
but it is an indefinite functional since it is cubically nonlinear.
A positive-definite functional that is conserved on smooth solutions is the momentum
\be
\mom =\sum_{k=1}^\infty  k \wphi^\ast(k,t) \wphi(k)
= \frac{1}{4\pi} \int_{\Ts} \left[\left|\partial_x\right|^{1/2} \phi(x,t)\right]^2\,dx.
\label{def:mom}
\ee
The \textit{a priori}
estimate provided by the conservation of $\mom$ (or the dissipation of $\mom$ for weak solutions)
is not sufficient to imply global existence, however, and 
we do not know of other conserved functionals of \eq{specee} for general $T$.

\subsection{Noncanonical variables}

It is convenient to introduce
noncanonical dependent variables
\[
u: \Ts \times I \to \Rl,\qquad \wu : \Ir\times I \to \Cx
\]
defined by
\[
u(x,t) = \left|\partial_x\right|^{1/2} \phi(x,t),\qquad \wu(k,t) = |k|^{1/2} \wphi(k,t).
\]
When rewritten in terms of $\wu$, equation \eq{specee} becomes
\be
\wu_t(k,t) + i k \sum_{n=-\infty}^\infty S(-k, k - n, n) \wu(k-n,t) \wu(n,t) = 0,
\label{ncspecee}
\ee
where
\be
S(k,m, n) = \frac{T(k,m, n)}{\left|k m n\right|^{1/2}}
\qquad \mbox{for $kmn \ne 0$}.
\label{defS}
\ee
We define $S(0,m,n) = S(k,0,n) = S(k,m,0) = 0$.
The corresponding spatial form of \eq{ncspecee} is then
\be
u_t + a(u,u)_x = 0,
\label{ncspatee}
\ee
where the bilinear form $a$ is given by
\be
\widehat{a(u,v)}(k,t) = \sum_{n=-\infty}^\infty
S(-k,k-n,n) \wu(k-n,t)\wv(n,t).
\label{defa}
\ee
The form $a$ in \eq{ncspatee} is related to $b$ in \eq{eeq} by
\be
b(\phi,\psi) = \left|\partial_x\right|^{1/2} a\left(\left|\partial_x\right|^{1/2}\phi,
\left|\partial_x\right|^{1/2}\psi\right).\label{ab}
\ee
It follows that equation \eq{ncspatee} is invariant under the rescalings
\[
t\mapsto \scaling t, \quad x\mapsto \scaling x,\quad
u \mapsto \scaling^{\nu-3/2} u.
\]

As we will show in Section~\ref{sec:propa}, under suitable assumptions on $S$ the
bilinear form $a(u,v)$ may be regarded as a kind of `nonlocal product'
of $u$ and $v$, so that \eq{ncspatee} is a nonlocal generalization of the
inviscid Burgers equation.
Our main result is that if the kernel $S(k,m,n)$ satisfies the bound in \eq{propS5} below,
where $\mu = \nu - 3/2 \ge 0$, and $s > \nu$, then the initial value problem for \eq{ncspatee}--\eq{defa} has a
unique solution for short times taking values in the Sobolev space
$H^s(\Ts)$. Equivalently, equation \eq{eeq}--\eq{defF} has a short-time $H^s$-solution
for $s > \nu + 1/2$.

\section{Examples}
\label{sec:ex}

In this section, we give a some examples of equations that belong to the class described in the previous section.
The simplest, and standard, example is the inviscid Burgers equation, but this class of equations
is surprising rich, and includes other interesting PDEs and integro-differential equations.
These examples also illustrate various different possibilities for local and global existence.

\subsection{Burgers equation}

Suppose that
\be
T(k,m,n) = |k m n|^{1/2}.
\label{burgker}
\ee
Then $S(k,m,n) = 1$ for $kmn\ne 0$, and
the noncanonical equation \eq{ncspatee} corresponding to \eq{burgker}
is the inviscid Burgers equation
\be
u_t+ \left(u^2\right)_x = 0,
\label{burgers}
\ee
as can be seen from \eq{defa} by use of the convolution theorem.
The kernel \eq{burgker} is homogeneous of degree $3/2$, corresponding to $n=d$ in
\eq{defnu}, and this equation describes bulk waves, such as weakly nonlinear
sound waves.

Theorem~\ref{th:locex} states local existence for
\eq{burgers} in $H^s$ with $s > 3/2$, which is the usual short-time existence result for
first-order quasilinear hyperbolic equations in one space dimension.
Nontrivial periodic solutions form singularities in finite time, but
the inviscid Burgers equation has global dissipative weak solutions that satisfy
appropriate entropy conditions. These solutions are not, however, compatible
with the Hamiltonian structure of \eq{burgers}.

\subsection{Hunter-Saxton equation}

Consider the kernel
\be
T(k, m, n) = \frac{1}{2}\left|k m n\right|^{1/2}\left(\frac{1}{ik} +\frac{1}{im} +  \frac{1}{in}\right)
\quad \mbox{for $kmn\ne 0$},
\label{hsker}
\ee
with $T(k,m,n) =0$ when $kmn=0$. Then
\[
S(k, m, n) = \frac{1}{2}\left(\frac{1}{ik} +\frac{1}{im} +  \frac{1}{in}\right)
\quad \mbox{for $kmn\ne 0$}.
\]
The corresponding noncanonical equation \eq{ncspatee} is
\be
u_t + \left(u \partial_x^{-1}u\right)_x - \frac{1}{2} \left(u^2 - \left\langle u^2\right\rangle\right)
 = 0,
\label{derivhseq}
\ee
where $\partial_x^{-1}$ is the anti-differentiation operator on zero-mean, periodic functions,
and the angular brackets denote the spatial average over a period. It follows that
$v = \partial_x^{-1}u$ satisfies the Hunter-Saxton equation \cite{saxton},
\be
v_t + vv_x - \frac{1}{2} \partial_x^{-1}\left(v_x^2
- \left\langle v_x^2\right\rangle\right) = 0.
\label{hseq}
\ee

The kernel \eq{hsker} has degree $\nu = 1/2$, corresponding to $n=d-2$ in \eq{defnu}. This is consistent with the
derivation of \eq{hseq} in \cite{saxton} for orientation waves in a massive director field. The waves propagate
in three-dimensional space, so $d=3$, and the mass parameter of the medium is set by the
moment of inertia per unit volume of the director field, which has dimension $M/L$, so $n=1$.

Theorem~\ref{th:locexneg} (with $\lambda=1$) implies the short-time existence of $H^s$-solutions of
\eq{derivhseq} for $s > 1/2$, and hence $H^s$-solutions of \eq{hseq}
for $s > 3/2$. This agrees with the short-time existence result of Yin \cite{yin}.
Nontrivial periodic solutions of the Hunter-Saxton equation form singularities in finite
time, but the equation has global weak solutions, including both dissipative
and energy-conserving weak solutions (\cite{bressan}, \cite{zheng}, \cite{zhang}).
The conservative weak solutions are compatible with the
Hamiltonian structure of the equation.

It is interesting to note that \eq{hseq} is completely integrable (\cite{beals}, \cite{zheng1}), and also
arises as the high-frequency limit of the integrable Camassa-Holm equation \cite{camassa},
\[
v_t + 2\kappa v_x + 3v v_x =
\left(v_t + vv_x\right)_{xx} - \frac{1}{2} \left(v_x^2\right)_x.
\]

\subsection{Surface wave equations}
\label{sec:surface}

As discussed in the introduction, weakly nonlinear asymptotics for `genuine' hyperbolic surface waves typically leads to
nonlocal scale-invariant equations. These equations inherit a Hamiltonian structure from the Hamiltonian
structure of the primitive physical equations (such as hyperelasticity or magneto\-hydrodynamics).

Perhaps the simplest kernel arising for surface waves is given by
\be
T(k,m,n) =
\frac{2|k m n|}{|k|+|m|+|n|},
\label{surfker}
\ee
leading to the spectral evolution equation
\be
\wphi_t(k,t) + i\sgn k \sum_{n=-\infty}^\infty
\frac{2|k| |k-n| |n|}{|k|+|k-n|+|n|} \wphi(k-n,t) \wphi(n,t) = 0.
\label{specsurfeq}
\ee
This equation was proposed by Hamilton \textit{et al}
(\cite{ham1}, \cite{ham2}) as a model equation for nonlinear Rayleigh waves, and it
describes the propagation of surface waves
on a tangential discontinuity in incompressible magnetohydrodynamics \cite{ali}.

The kernel in \eq{surfker} is homogeneous of degree $\nu=2$. This degree is consistent with the result of
dimensional analysis in \eq{defnu}, since the dimension $d$ of the wavenumber space
of surface waves that propagate along the boundary is one
less that the spatial dimension $n$ of the medium (whose density sets the mass parameter of the
wave motion).

The spatial form of \eq{specsurfeq} may be written explicitly as (\cite{parker}, \cite{ham2})
\be
\phi_t  +\frac{1}{2}
\hilbert\left[\psi^2\right]_{xx}
+ \psi \phi_{xx}
 = 0,
 \label{modelee}
 \ee
where $\hilbert$ is the Hilbert transform and $\psi = \hilbert[\phi]$. This is equation \eq{surf_ex}.

Equation \eq{modelee} corresponds to \eq{eeq} with
\beann
b(\phi,\phi) &=& \frac{1}{2}\left(\psi^2\right)_{xx} - \hilb\left[\psi \phi_{xx}\right]\\
&=& - \frac{\delta}{\delta\phi} \int \psi\psi_x\phi_x\, dx.
\eeann
The appearance of second-order spatial derivatives in the expression for $b$ is misleading
since some cancelation occurs. It follows from \eq{ab} and Proposition~\ref{prop:a} (with $\mu=1/2$) that
$b : H^s \times H^s \to H^{s-1}$
involves the `loss' of only one derivative when $s > 3/2$.

Theorem~\ref{th:locex} implies that the noncanonical equation
associated with \eq{modelee} has short-time $H^s$-solutions when $s > 2$, so
\eq{modelee} has $H^s$-solutions when $s > 5/2$.

Numerical computations (\cite{ali}, \cite{ham2}, \cite{parker3})
show that solutions of \eq{modelee} form singularities
in which the derivative $\phi_x$ blows up, but $\phi$ appears to remain continuous.
The global existence of appropriate weak solutions is an open question.
As far as we have been able to determine, equation \eq{modelee} is not completely integrable.

Similar equations describe the evolution of weakly nonlinear Rayleigh waves
on an elastic half-space.
In that case, $\phi(x,t)$ corresponds to the vertical displacement of the boundary of the half-space.
For weakly nonlinear Rayleigh waves on an isotropic elastic half-space, $\wphi(k,t)$
satisfies an equation of the form \eq{specee} with the kernel
(see \cite{kal1}, \cite{kal2}, \cite{lardner}, \cite{mayer}, \cite{parker1},
\cite{zabolotskaya}, for example)
\bea
T(k,n,m) &=& \frac{\alpha|k| |m| |n|}{|k| +r|m| + r|n|}
+  \frac{\alpha|k| |m| |n|}{r|k| +|m| + r|n|}
\nonumber\\
&&+\frac{\alpha|k| |m| |n|}{r|k| +r|m| + |n|}
+\frac{\beta|k| |m| |n|}{r|k| +|m| + |n|}
\nonumber\\
&&+\frac{\beta|k| |m| |n|}{|n| +r|m| + |n|}
+ \frac{\beta|k| |m| |n|}{|k| +|m| + r|n|}
\nonumber\\
&&+\frac{\gamma|k| |m| |n|}{|k| +|m| + |n|}.
\label{rayleighker}
\eea
Here, $0 < r = c_t/c_l < 1$ is the ratio of the elastic
solid's transverse and longitudinal wave speeds, $c_t$ and $c_l$ respectively,
and $\alpha$, $\beta$, $\gamma$ are nonlinear elastic constants.
The corresponding equation for $\phi(x,t)$ has short-time $H^s$-solutions for $s > 5/2$.

More complicated, but qualitatively similar, kernels arise for nonlinear Rayleigh
waves on a non-isotropic elastic half-space (\cite{ham3}, \cite{lardner1}).
Kernels with a different degree of homogeneity ($\nu=5/2$) arise in the weakly nonlinear description of
elastic edge waves (\cite{kiselev}, \cite{parker2}).

\subsection{Compacton equations}
\label{sec:compacton}

The kernel
\be
T(k, m, n) =  i \left|k m n\right|^{1/2} k m n
\label{compactonker}
\ee
leads to a noncanonical equation \eq{ncspatee} that is a local PDE,
\[
u_t + \left(u_x^2\right)_{xx} = 0.
\]
Differentiating this equation with respect to $x$ and setting $u_x=v$, we get
\be
v_t + \left(v^2\right)_{xxx} = 0.
\label{hfkmneq}
\ee
This equation is the high-frequency limit of the $K(2,2)$-equation
introduced by Rosenau and Hyman \cite{rosenau}:
\be
v_t + \left(v^2\right)_{x} + \left(v^2\right)_{xxx} = 0.
\label{compacton}
\ee
The $K(2,2)$-equation has `compacton' solutions
(that is, compactly supported traveling waves) given by
\[
v(x,t) = \frac{4 c}{3} \cos^2\left(\frac{x-ct}{4}\right) \qquad\mbox{if $|x-ct| \le 2\pi$},
\]
and $v(x,t) = 0$ if $|x-ct| >  2\pi$.
Some numerical solutions of \eq{compacton} are given in \cite{frutos}, \cite{levy}.

Although \eq{hfkmneq} does not have compacton solutions, it has travelling wave solutions
that describe the local non-smooth behavior of the compacton solutions of \eq{compacton}
near the `edges' of the compacton; for example:
\[
v(x,t) =
\frac{c}{12}(x-ct)^2\qquad \mbox{if $x\le ct$},
\]
and $v(x,t) = 0$ if $x > ct$.
The relationship between \eq{hfkmneq} and \eq{compacton} is thus similar to the one between the Hunter-Saxton
and Camassa-Holm equations, including the local description of the non-smooth behavior of the `peakon'
solutions of the Camassa-Holm equation by solutions of the Hunter-Saxton equation.
As far as we know, however, neither \eq{hfkmneq} nor \eq{compacton} is completely integrable.

An expansion of the spatial derivatives in \eq{hfkmneq},
\[
v_t + 2 v v_{xxx} + 6 v_x v_{xx} = 0,
\]
leads to a backward diffusion term when $v_x > 0$, which suggests that the
equation may not be well-posed. The short-time existence proof given here
does not apply to \eq{hfkmneq} since its kernel does not satisfy \eq{estimate}
below, and this condition seems to be necessary to obtain $L^2$-Sobolev estimates.
It is conceivable that other estimates could lead to short-time existence
for \eq{hfkmneq} or \eq{compacton}, but no such results appear to be
known (although it is possible that a Cauchy-Kowalevski type argument could be
used for analytic initial data \cite{costin}).

More generally, the kernel
\[
T(k, m, n) = \left|k m n\right|^{1/2} \left(ik m n\right)^{p}
\]
of degree $\nu = 3p + 3/2$, has the corresponding noncanonical equation
\[
u_t + \partial_{x}^{p+1} \left(\partial_x^p u\right)^2 = 0,
\]
which may be regarded as a higher-order generalization of the inviscid Burgers equation,
to which it reduces when $p=0$. Our short-time existence result does not apply to this equation when
$p > 0$.


\section{Local existence of smooth solutions}
\setcounter{equation}{0}
\label{sec:locex}

We consider an initial value problem for the noncanonical form
of the equation described in Section~\ref{sec:ee}:
\bea
&&u_t + a(u,u)_x = 0,
\label{ncspatee1}
\\
&&u(x,0) = f(x).
\nonumber
\eea
Here, $u : \Ts \times I \to \Rl$, with $I$ an open time interval containing the origin,
$f : \Ts \to \Rl$ is given initial data, and the bilinear form $a$ is defined by
\be
\widehat{a(u,v)}(k,t) = \sum_{n=-\infty}^\infty
S(-k,k-n,n) \wu(k-n,t)\wv(n,t).
\label{defa1}
\ee
We assume that the kernel $S :\Ir^3 \to \Cx$ satisfies
\bea
&&S(k,m,n) = 0\quad \mbox{if $kmn=0$},
\label{propS1}\\
&&S(k, m , n)  =S^\ast(-k, -m , -n),
\label{propS1a}\\
&&S(k, m , n) = S(k, n , m),
\label{propS2}\\
&&S(k, m , n) = S(m , n, k),
\label{propS3}\\
&&S(\scaling k, \scaling m , \scaling n) = \scaling^\mu S(k, m , n)\qquad \mbox{for $\scaling >0$},
\label{propS4}
\eea
These properties follow from those of $T$ in
\eq{propT} and the definition of $S$ in \eq{defS}, with
\[
\mu = \nu - \frac{3}{2}.
\]
We will focus our attention on kernels $S$ with nonnegative degree $\mu \ge 0$.
We consider kernels of negative degree in Section~\ref{sec:neg}.

In order to develop a local existence theory, we assume
that $S$ satisfies the following condition for
some constant $C$ and all $k,m,n\in \Ir$:
\be
|S(k, m, n)|
\le C\min\left\{ |k|^{\mu}, |m|^{\mu}, |n|^{\mu}\right\};
\label{propS5}
\ee
equivalently, $T$ satisfies
\be
{\left|T(k, m, n)\right|}
\le C \left|k m n\right|^{1/2}\min\left\{ |k|^{\nu-3/2}, |m|^{\nu-3/2}, |n|^{\nu-3/2}\right\}.
\label{estimate}
\ee

For example, the kernel \eq{burgker} for the inviscid Burgers equation satisfies \eq{estimate}
with $\nu=3/2$ and $C=1$. The surface-wave kernel \eq{surfker} satisfies \eq{estimate} with
$\nu=2$ and $C=1$, since
\[
\frac{2|k m n|^{1/2}}{|k|+|m|+|n|} \le\min\left\{ |k|^{1/2}, |m|^{1/2}, |n|^{1/2}\right\},
\]
and the Rayleigh-wave kernel \eq{rayleighker} satisfies \eq{estimate} with $\nu=2$ and
\[
C = \frac{1}{2}\left(\frac{3|\alpha|}{r} +  \frac{3|\beta|}{r} + |\gamma|\right).
\]
By contrast, the compacton-type kernel \eq{compactonker} does not satisfy \eq{estimate}.

\subsection{Notation}

We denote by $\dot H^s(\Ts)$, or $\dot H^s$ for short, the Hilbert space
of periodic functions with zero mean and square integrable
derivatives of the order $s$,
\[
\dot H^s(\Ts) = \left\{f : \Ts \to \Rl \mid \sum_{k=-\infty}^\infty |k|^{2s}
\left|\wf(k)\right|^2 < \infty,
\quad \wf(0) = 0\right\}.
\]
As inner product and norm on $\dot H^s$, we use
\beann
\langle f, g\rangle_s
&=& \sum_{k=-\infty}^\infty |k|^{2s} \wf(k)\wg(-k),\\
\|f\|_s
&=& \left(\sum_{k=-\infty}^\infty |k|^{2s} \left|\wf(k)\right|^2\right)^{1/2}.
\eeann
In particular,
\[
\|f\|_0 = \left(\frac{1}{2\pi} \int_{\Ts} |f(x)|^2\, dx\right)^{1/2}.
\]
For $\wf : \Ir\to \Cx$, we also define
\[
\left\|\wf\right\|_{\ell^p} = \left(\sum_{k=-\infty}^\infty \left|\wf(k)\right|^p\right)^{1/p}
\]
so that
\be
\|f\|_s = \Bigl\| |k|^{s} \wf\Bigr\|_{\ell^2}.
\label{sobineq1}
\ee

\subsection{The bilinear form}
\label{sec:propa}

First, we obtain some properties of the bilinear form $a$ given by \eq{defa1}.
In the simplest case when $S(k,m,n) =1$ for $kmn\ne0$, corresponding to $\mu=0$, the
form $a$ is a multiplication operator on $\dot H^s$ for $s > 1/2$:
\be
a(u,v) = uv - \frac{1}{2\pi}\int_{\Ts} uv \, dx.
\label{pointprod}
\ee
As the following proposition shows, the form $a$
has similar properties to a multiplication operator
for any kernel $S$ that satisfies \eq{propS5}, although we need some additional
smoothness to ensure that $\dot H^s$ is an algebra under the product defined by $a$.

\begin{proposition}
\label{prop:a}
Let $a(u,v)$ be defined by \eq{defa1}, where $S: \Ir^3 \to \Cx$ satisfies \eq{propS1}--\eq{propS5}.
Suppose that $s > \mu + 1/2$.

\smallskip\noindent
(a) Then $a : \dot H^s\times \dot H^s \to \dot H^{s}$
is a bounded symmetric bilinear form.

\smallskip\noindent
(b) For all $u,v\in \dot H^{s+1}$,
\be
a(u,v)_x = a(u_x,v) + a(u,v_x).
\label{deriv}
\ee

\smallskip\noindent
(c) For all $u, v, w \in \dot H^{s}$,
\be
\frac{1}{2\pi}\int_{\Ts} u a(v,w)\,dx = \frac{1}{2\pi}\int_{\Ts} v a(w,u)\,dx,
\label{cyclicint}
\ee
and there is a constant $M_{s-\mu} > 0$ such that
\be
\left|\frac{1}{2\pi}\int_{\Ts} u a(v,w)\,dx\right| \le C M_{s-\mu} \|u\|_s \|v\|_0 \|w\|_0.
\label{cyclicbound}
\ee
\end{proposition}

\medskip\noindent
\textbf{Proof.} Using \eq{defa1}, we get
\bea
\left\|a(u,v)\right\|_s^2 &=& \sum_{k} |k|^{2s} \left|\widehat{a(u,v)}(k)\right|^2
\nonumber\\
&=& \sum_{k} \left| |k|^{s} \sum_{n} S(-k,k-n,n) \wu(k-n)\wv(n)\right|^2.
\label{randomeq}
\eea
For $s > 0$ and $m,n \in \Rl$, we have
\[
\left|m+n\right|^{s} \le B_s\left(|m|^s + |n|^{s}\right),
\]
where
\[
B_s = \left\{\begin{array}{ll}
1 & \mbox{if $0 < s \le 1$},\\
2^{s-1} & \mbox{if $s \ge 1$}.
\end{array}\right.
\]
Using this inequality, with $m = k-n$, and \eq{propS5} in \eq{randomeq}, we get
\beann
\left\|a(u,v)\right\|_s^2&\le& B_s^2 \sum_{k} \left|\sum_{n} \left(|k-n|^{s} + |n|^{s}\right) \left|S(-k,k-n,n) \wu(k-n)\wv(n)\right|\right|^2
\\
&\le& B_s^2 C^2 \sum_{k} \left|\sum_{n}\left(|k-n|^{s}|n|^{\mu} + |k-n|^{\mu}|n|^{s}
\right)\left|\wu(k-n)\wv(n)\right|\right|^2
\\
&\le& B_s^2 C^2\Bigl\|\left(|k|^s |\wu|\right)\ast \left(|k|^\mu |\wv|\right)
+ \left(|k|^\mu |\wu|\right)\ast \left(|k|^s |\wv|\right)\Bigr\|_{\ell^2}^2.
\eeann
Here, we denote the convolution of $\wf, \wg : \Ir \to \Cx$ by
\[
\wf\ast\wg(k) = \sum_{n=-\infty}^{\infty} \wf(k-n)\wg(n).
\]
It follows that
\[
\left\|a(u,v)\right\|_s \le B_s C\left\{\Bigl\|\left(|k|^s |\wu|\right)\ast \left(|k|^\mu |\wv|\right)\Bigr\|_{\ell^2}
+ \Bigl\|\left(|k|^\mu |\wu|\right)\ast \left(|k|^s |\wv|\right)\Bigr\|_{\ell^2}\right\}.
\]
Using Young's inequality,
\be
\left\|\wf\ast\wg\right\|_{\ell^2} \le \Bigl\|\wf\Bigr\|_{\ell^1}\left\|\wg\right\|_{\ell^2},
\label{youngineq}
\ee
we get that
\be
\left\|a(u,v)\right\|_s \le B_s C  \left\{ \|u\|_s \Bigl\| |k|^\mu |\wv|\Bigr\|_{\ell^1}
+  \Bigl\| |k|^\mu |\wu|\Bigr\|_{\ell^1} \|v\|_s\right\}.
\label{abound}
\ee
If $s > 1/2$, we have for all $f\in \dot H^s$ that
\be
\left\|\wf\right\|_{\ell^1} \le M_{s} \left\|f\right\|_s,
\label{sobineq}
\ee
where
\be
M_s = \left(2\sum_{n=1}^\infty \frac{1}{n^{2s}}\right)^{1/2}.
\label{defms}
\ee
Hence, if $s > \mu + 1/2$, inequalities \eq{abound} and \eq{sobineq} imply that
\[
\|a(u,v)\|_s \le 2 B_s C M_{s-\mu} \|u\|_s \|v\|_s.
\]
Thus, $a$ is bounded on $\dot H^s$. The symmetry and bilinearity of $a$ is obvious.

To prove \eq{deriv}, we observe that
\beann
\widehat{a(u,v)_x}(k) &=& ik  \widehat{ a(u,v)}(k) \\
&=& ik  \sum_{n} S(-k,k-n,n) \wu(k-n)\wv(n)\\
&=&\sum_{n} i(k-n)S(-k,k-n,n) \wu(k-n)\wv(n)\\
&&\qquad + \sum_{n} in S(-k,k-n,n) \wu(k-n) \wv(n)\\
&=& \widehat{ a(u_x,v)}(k) + \widehat{ a(u,v_x)}(k).
\eeann

Equation \eq{cyclicint} follows from the fact that
\beann
\frac{1}{2\pi} \int_{\Ts} u a(v,w)\,dx &=& \sum_k \wu(-k) \widehat{a(v,w)}(k) \\
&=& \sum_{k,n} S(-k, k-n, n) \wu(-k) \wv(k-n) \ww(n)
\eeann
is symmetric under cyclic permutations of ${u,v,w}$. Moreover, using \eq{propS5},
the Cauchy-Schwartz inequality, Young's inequality, and \eq{sobineq} in this
equation, we get
for $s > \mu + 1/2$ that
\beann
\left|\frac{1}{2\pi} \int_{\Ts} u a(v,w)\,dx \right| &\le& C\sum_{k,n} |k|^\mu \left|\wu(-k) \wv(k-n) \ww(n)\right|\\
&\le& C\Bigl\|\left(|k|^\mu |\wu|\right)\ast \left(|\wv|\right)\Bigr\|_{\ell^2} \left\|\ww\right\|_{\ell^2}\\
&\le& C\Bigl\||k|^\mu \wu\Bigr\|_{\ell^1} \left\|\wv\right\|_{\ell^2}\left\|\ww\right\|_{\ell^2}\\
&\le& C M_{s-\mu}\left\|u\right\|_{s} \left\|v\right\|_{0}\left\|w\right\|_{0}.\qquad\Box
\eeann


\bigskip
The form $a$ satisfies other Sobolev inequalities. In Proposition~\ref{prop:a}, we prove only the ones that
we need to use below.

\subsection{Local existence}

We define
\be
K_s = C C_s M_{s-\mu-1},
\label{defks}
\ee
where $C$ is the constant in \eq{propS5}, $C_s$ is the constant in \eq{ineqcs}, and $M_s$
is given by \eq{defms}.

\begin{theorem}
\label{th:locex}
Suppose that $S : \Ir^3 \to \Cx$ satisfies \eq{propS1}--\eq{propS4} and \eq{propS5}
with $\mu \ge 0$, the form $a$ is given by \eq{defa1}, and $s > \mu + 3/2$.
Then, for any $f\in \dot H^s(\Ts)$, the initial-value problem \eq{ncspatee1}
has a unique local solution
\be
u \in C\left(I; \dot H^s(\Ts)\right) \cap C^1\left(I;\dot H^{s-1}(\Ts)\right)
\label{propu}
\ee
defined on the time interval $I = (-t_\ast,t_\ast)$, where
\be
t_\ast =  \frac{1}{K_s \left\|f\right\|_s}.
\label{deftstar}
\ee
\end{theorem}

\medskip\noindent
\textbf{Proof.}
Let $P^N : \dot H^s(\Ts) \to \dot H^s(\Ts)$ denote the orthogonal projection
\be
P^N \left(\sum_{k=-\infty}^\infty \wf(k) e^{ikx}\right) = \sum_{k=-N}^N \wf(k) e^{ikx}.
\label{proj}
\ee
We introduce the Galerkin approximations

\[
u^N(x,t) = \sum_{k=-N}^N \wu^N(k,t) e^{ikx}
\]
that satisfy the system of ODEs
\bea
&&u^N_t + P^N a\left(u^N,u^N\right)_x = 0,
\label{galerkin}\\
&&u^N(x,0) = P^N f(x).
\nonumber
\eea
The main estimate for the local existence proof is given in the following
Proposition.

\begin{proposition}
\label{prop:estimate}
Suppose $u^N(x,t)$ satisfies \eq{galerkin} and $s > \mu + 3/2$. Then
\be
\left|\frac{d}{dt} \left\|u^N\right\|_s\right| \le
K_s\left\|u^N\right\|_s^2,
\label{nlineq}
\ee
where $K_s$ is defined in \eq{defks}.
\end{proposition}

\medskip\noindent
\textbf{Proof.}
Dropping the superscript $N$ to simplify the notation, we compute from \eq{galerkin}, \eq{proj}, and \eq{defa1}
that
\bea
&&\frac{d}{dt} \sum_{k} |k|^{2s} \wu(k,t)\wu(-k,t)
\nonumber\\
&&+ 2i \sum_{k,n} k|k|^{2s} S(-k,k-n,n)  \wu(- k,t)\wu(k-n,t)\wu(n,t) = 0.
\label{basiceq2}
\eea
Using the cyclic symmetry of $S(k,m,n)$, we can rewrite this equation as
\beann
&&\frac{d}{dt} \sum_{k} |k|^{2s} \left|\wu(k,t)\right|^2
\\
&&\qquad+ \frac{2i}{3} \sum_{k,n} \left\{
k|k|^{2s} - (k-n)|k-n|^{2s} - n|n|^{2s}\right\}
\\
&&\qquad\qquad\quad S(-k,k-n,n)\wu(- k,t)\wu(k-n,t)\wu(n,t) = 0.
\eeann
Using the inequality \eq{ineqcs}, we find that
\bea
&&\left|\frac{d}{dt} \sum_{k} |k|^{2s}  \left|\wu(k,t)\right|^2\right|
\nonumber\\
&&\qquad\qquad \le
\frac{2 C_s}{3} \sum_{k,n} \left\{|k|^{s}|k-n|^{s}|n|+ |k-n|^{s}|n|^{s}|k|+ |n|^{s}|k|^{s}|k-n|\right\}
\nonumber\\
&&\qquad\qquad\qquad\qquad \left|S(-k,k-n,n)\right|\,\left|\wu(- k,t)\wu(k-n,t)\wu(n,t)\right|.
\label{basiceq}
\eea
Using \eq{propS5} in the right hand side of \eq{basiceq}, and applying the appropriate
bound on each term, we get
\beann
&&\left|\frac{d}{dt} \sum_{k} |k|^{2s}  \left|\wu(k,t)\right|^2\right|
\\
&&\qquad\qquad \le\frac{2 C C_s}{3} \sum_{k,n} \left\{|k|^{s}|k-n|^{s}|n|^{\mu + 1}
+ |k-n|^{s}|n|^{s}|k|^{\mu + 1}\right.
\\
&&\qquad\qquad \qquad\qquad\left.+ |n|^{s}|k|^{s}|k-n|^{\mu + 1}
\right\} \left|\wu(- k,t)\wu(k-n,t)\wu(n,t)\right|
\\
&&\qquad\qquad \le {2 C C_s} \sum_{k,n} |k|^{s}|k-n|^{s}|n|^{\mu + 1} \left|\wu(- k,t)\wu(k-n,t)\wu(n,t)\right|.
\eeann
Using the Cauchy-Schwartz inequality and Young's inequality \eq{youngineq} on the
right hand side of this inequality, we get
\beann
\left|\frac{d}{dt} \Bigl\||k|^{s} \wu\Bigr\|_{\ell^2}^2\right|&\le&
{2 C C_s} \Bigl\| |k|^{s}\wu\Bigr\|_{\ell^2}
\Bigl\|(|k|^{s}|\wu|)\ast(|k|^{\mu + 1} |\wu|)\Bigr\|_{\ell^2},
\\
&\le&
{2 C C_s} \Bigl\| |k|^{s}\wu\Bigr\|_{\ell^2}^2 \Bigl\||k|^{\mu + 1} \wu\Bigr\|_{\ell^1}.
\eeann
It then follows from \eq{sobineq1} and \eq{sobineq} that if $s > \mu + 3/2$, then
\[
\left|\frac{d}{dt} \| u\|_s^2\right|\le
{2 C C_s M_{s-\mu-1}} \| u\|_s^3.
\]
Simplifying this equation and using \eq{defks} we get \eq{nlineq}. $\Box$

\bigskip
Using Gronwall's inequality, we deduce from \eq{nlineq} that
\[
\left\|u^N\right\|_s(t) \le \frac{\left\|f\right\|_s}{1 - K_s \left\|f\right\|_s |t|},
\]
so the Galerkin approximations exist in $|t| < t_\ast$, where
$t_\ast$ is given by \eq{deftstar}.

It follows from Proposition~\ref{prop:a} and standard arguments
(for example, \cite{taylor}, \S 16.1, or \cite{kato})
that we can extract a subsequence $(u^N)$ of Galerkin approximations such that
\be
u^N \to u \in C_w\left(J; \dot H^s\right)
\label{provreg}
\ee
on any closed interval $J \subset{I}$, where $u$ is a distributional solution of \eq{ncspatee1}.
(Here, $C_w$ denotes the space of weakly continuous functions.)
Proposition~\ref{prop:estimate} and Gronwall's inequality
imply that $t\mapsto \|u(t)\|_s$ is continuous. It follows that
$u$ is strongly $\dot H^s$-continuous, and thus satisfies \eq{propu}
since $a$ is strongly continuous on $\dot H^s$.

To prove uniqueness, we suppose that $u$, $v$ are two solutions of \eq{ncspatee1} with the regularity stated
in \eq{propu}. Then, using integration by parts and the properties
of $a$ given in Proposition~\ref{prop:a}, we compute that
\beann
\frac{d}{dt}  \left\|u-v\right\|_0^2 &=& \frac{1}{2\pi}\frac{d}{dt} \int_{\Ts} \left(u-v\right)^2\,dx\\
&=& - \frac{1}{\pi}\int_{\Ts} \left(u-v\right) \left[a(u,u)_x-a(v,v)_x\right]\,dx\\
&=& - \frac{1}{\pi}\int_{\Ts} \left(u-v\right) a(u+v,u-v)_x\,dx\\
&=& \frac{1}{\pi}\int_{\Ts} \left(u-v\right)_x a(u+v,u-v) \,dx\\
&=& \frac{1}{\pi}\int_{\Ts} \left(u+v\right) a\left(u-v,u_x - v_x\right) \,dx\\
&=& \frac{1}{2\pi}\int_{\Ts} \left(u+v\right) a\left(u-v,u - v\right)_x \,dx\\
&=& -\frac{1}{2\pi}\int_{\Ts} \left(u+v\right)_x a\left(u-v,u - v\right) \,dx.
\eeann
Since $s > \mu + 3/2$, it follows from this equation and \eq{cyclicbound} that
\[
\left|\frac{d}{dt} \left\|u-v\right\|_0^2 \right| \le C M_{s-\mu-1} \|u + v\|_s   \left\|u-v\right\|_0^2.
\]
Hence, $u=v$ by Gronwall's inequality, and the solution is unique. $\Box$

\subsection{Negative degree kernels}
\label{sec:neg}

In this section, we consider kernels $S$ of negative degree
$\mu=-\lambda < 0$ that satisfy the following condition
for some constant $C$ and all $kmn\ne 0$:
\be
|S(k, m, n)|
\le \frac{C}{\min\left\{|k|^{\lambda}, |m|^{\lambda}, |n|^{\lambda}\right\}}.
\label{propS5neg}
\ee
This corresponds to the assumption that $T$ satisfies the condition
\be
{\left|T(k, m, n)\right|}
\le \frac{C \left|k m n\right|^{1/2}}{\min\left\{|k|^{3/2-\nu}, |m|^{3/2-\nu}, |n|^{3/2-\nu}\right\}}.
\label{negestimate}
\ee
where $\nu = 3/2-\lambda < 3/2$.
In particular, \eq{negestimate} is satisfied by the Hunter-Saxton kernel \eq{hsker}
with $\nu=1/2$ and $C= 3/2$.

Let
\[
K_{s,\lambda} = C C_{s,\lambda} M_{s+\lambda-1},
\]
where $C$ is the constant in \eq{propS5neg}, $C_{s,\lambda}$ is the constant in \eq{ineqds}, and $M_s$
is given by \eq{defms}.

\begin{theorem}
\label{th:locexneg}
Suppose that $S : \Ir^3 \to \Cx$ satisfies \eq{propS1}--\eq{propS4} and \eq{propS5neg}
with $\lambda > 0$, the form $a$ is given by \eq{defa1}, and
\[
s > \max\left\{\frac{3}{2} - \lambda,\frac{1}{2}\right\}.
\]
Then, for any $f\in \dot H^s(\Ts)$, the initial-value problem \eq{ncspatee1}
has a unique local solution
\[
u \in C\left(I; \dot H^s(\Ts)\right) \cap C^1\left(I;\dot H^{s-1}(\Ts)\right)
\]
defined on the time interval $I = (-t_\ast,t_\ast)$, where
\[
t_\ast =  \frac{1}{K_{s,\lambda} \left\|f\right\|_s}.
\]
\end{theorem}

\medskip\noindent
\textbf{Proof.} The proof is essentially the same as before. For $s > 1/2$,
we get that $a: \dot H^s\times \dot H^s \to \dot H^s$ is a continuous bilinear form.
From \eq{basiceq2} and \eq{propS5neg}, the Galerkin approximations satisfy
\beann
&&\left|\frac{d}{dt} \sum_{k} |k|^{2s}  \left|\wu(k,t)\right|^2\right|
\nonumber\\
&&\quad\le
\frac{2}{3} \sum_{k,n} \left| k|k|^{2s} - (k-n)|k-n|^{2s} - n|n|^{2s}\right|
\nonumber\\
&&\qquad\qquad\quad\left|S(-k,k-n,n)\right| \left|\wu(- k,t)\wu(k-n,t)\wu(n,t)\right|
\nonumber\\
&&\quad\le
\frac{2C}{3} \sum_{k,n} \frac{\left|k|k|^{2s} - (k-n)|k-n|^{2s} - n|n|^{2s}\right|}{ \min\left\{|k|^\lambda,
|k-n|^\lambda,|n|^{\lambda}\right\}}\left|\wu(- k,t)\wu(k-n,t)\wu(n,t)\right|.
\eeann
Using \eq{ineqds} in this equation, and summing over cyclic permutations,
we find that
\beann
\left|\frac{d}{dt} \sum_{k} |k|^{2s}  \left|\wu(k,t)\right|^2\right|
&\le&
{2 C C_{s,\lambda}} \sum_{k,n} |k|^{s}|k-n|^{s}|n|^{1-\lambda}
\\
&&\qquad\qquad\quad \left|\wu(- k,t)\wu(k-n,t)\wu(n,t)\right|
\\
&\le& {2 C C_{s,\lambda}} \Bigl\||k|^{1-\lambda} \wu\Bigr\|_{\ell^1} \Bigl\| |k|^{s}\wu\Bigr\|_{\ell^2}^2.
\eeann
Hence, using \eq{sobineq}, with $s > 3/2-\lambda$, and simplifying the result, we get
\[
\left|\frac{d}{dt} \left\|u\right\|_s \right| \le {C C_{s,\lambda} M_{s+\lambda-1}}\left\|u\right\|_s^2 ,
\]
and the local existence result follows. $\Box$

\section{An Inequality}
\label{sec:ineq}
\setcounter{equation}{0}

For completeness, we prove an inequality used in our local existence proof.

\begin{proposition}
\label{prop:ineq}
For $s > 0$ and $\lambda \ge 0$, there exists a constant $C_{s,\lambda}$ such that
for all nonzero $k,m,n \in \Rl$ with $k+m+n=0$,
\bea
&&\frac{\left|k|k|^{2s} + m|m|^{2s} + n|n|^{2s}\right|}{ \min\left\{|k|^\lambda,
|m|^\lambda,|n|^{\lambda}\right\}}
\nonumber\\
&&\quad
\le C_{s,\lambda} \left\{|k|^{s}|m|^{s}|n|^{1-\lambda}
+ |n|^{s}|k|^{s}|m|^{1-\lambda} + |m|^{s}|n|^{s}|k|^{1-\lambda}\right\}.
\label{ineqds}
\eea
\end{proposition}

\medskip\noindent
\textbf{Proof.}
We can assume without loss of generality that $0 < m \le n$ and $k = -(m+n) < 0$.
Writing $x = m/n$, we have
\[
\frac{\left|k|k|^{2s} + m|m|^{2s} + n|n|^{2s}\right|}{|m|^{\lambda} \left\{|k|^{s}|m|^{s}|n|^{1-\lambda}
+ |n|^{s}|k|^{s}|m|^{1-\lambda} + |m|^{s}|n|^{s}|k|^{1-\lambda}\right\}} = f_{s,\lambda}(x),
\]
where $f_{s,\lambda} : (0,1] \to (0,\infty)$ is given by
\[
f_{s,\lambda}(x) = \frac{(x+1)^{2s+1} - x^{2s+1} - 1}{(x+1)^s \left(x^{s+\lambda} + x\right)
+ (x+1)^{1-\lambda}x^{s+\lambda}}.
\]
This function is continuous on $(0,1]$ and has a finite limit as $x\to 0^+$, with
\[
f_{s,\lambda}\left(0^+\right) = \left\{\begin{array}{ll}
0 & \mbox{if $0 < s + \lambda < 1$},
\\
(2s+1)/3 & \mbox{if $ s + \lambda = 1$},
\\
2s+1 & \mbox{if $s + \lambda > 1 $}.
\end{array}\right.
\]
Hence $f_{s,\lambda}$ is bounded on $(0,1]$
by $C_{s,\lambda}$, say, which proves \eq{ineqds}. $\Box$

\bigskip
In particular, if $\lambda=0$, we write $f_s= f_{s,0}$ and $C_s = C_{s,0}$, so that
\bea
&&\left|k|k|^{2s} + m|m|^{2s} + n|n|^{2s}\right|
\nonumber\\
&&\qquad
\le C_s \left\{|k|^{s}|m|^{s}|n| + |n|^{s}|k|^{s}|m| + |m|^{s}|n|^{s}|k|\right\}
\label{ineqcs}
\eea
for all $k,m,n \in \Rl$ with $k+m+n=0$.
We have $f_1(x) =1$, $f_3(x) = 7$, and $f_s(1) = 2^s-1$.
Numerical plots show that $f_s$ is monotone decreasing on $(0,1]$ for
$1< s \le 3$ and monotone increasing for $s \ge 3$, in which case
\be
C_s = \left\{\begin{array}{ll}
2s+1 & \mbox{if $1 < s \le 3$},
\\
2^s-1 & \mbox{if $s \ge 3$}.
\end{array}\right.
\label{defcs}
\ee

\end{document}